\documentstyle[10pt,amsfonts,german]{article}
\newtheorem{theorem}{Theorem}[section]
\newtheorem{satz}[theorem]{Theorem}  
\newtheorem{lem}[theorem]{Lemma}  
\newtheorem{cor}[theorem]{Corollary}  
\newtheorem{defi}[theorem]{Definition}  
\newcommand{\titel}[1]{\begin{center} {\LARGE\bf\sc\ #1}
\end{center}\title{#1}\vspace{1cm}}
\newcommand{\autor}[1]{\begin{center}{\sc #1}\end{center} \author{#1}}
\newcommand{\adresse}[1]{ \begin{center} {\it #1}
\end{center}\vspace{0,3cm}}
\newcommand{\abstrakt}[4]{\mbox{} \hrulefill \mbox{} \begin{quotation}\noindent 
{\footnotesize\bf Abstract: }{\footnotesize #1 }\\[0,5cm] 
 {\footnotesize\it Subj. Class:} {\footnotesize #2}\\
  {\footnotesize\it 1991 MSC:} {\footnotesize #3} \\ {\footnotesize\it Keywords:}
 {\footnotesize #4}
 \end{quotation} 
 \vspace{0,1cm}
\mbox{} \hrulefill \mbox{}\vspace{1cm}}
\newcommand{\map}[3]{#1:#2\longrightarrow #3}  
\newcommand{\twist}{{\cal T}}  
\newcommand{\dirac}{{\cal D}}
\newcommand{\R}{{\cal R}}
\newcommand{\krumm}{{\cal R}}

\newcommand{\r}{ {\mathbb{R} } }
\newcommand{\com}{ { \mathbb{C}} }
\newcommand{\la}[1]{ {\mathfrak{ #1 }} }
\newcommand{\dual}{(\r^n)^*}

\newcommand{\beweis}{ {\sc Proof:} }
\newcommand{\bemerkung}{ {\bf Remark:} }

\newcommand{\schluss}{\nopagebreak  \hfill $\Box$ \\[0,2cm] }  
\newcommand{\gleifeldnon}[2]{\begin{eqnarray*}\lefteqn{#1}\\#2 \end{eqnarray*}}

\newcommand{\glungen}[1]{\begin{eqnarray*}#1 \end{eqnarray*}}
\newcommand{\lab}[1]{\label{#1}}
\newcommand{\glei}[2]{\begin{equation}\lab{#1} #2 \end{equation}}
\newcommand{\dichte}{ {\cal L}}
\newcommand{\sa}[1]{\begin{theorem}#1\end{theorem}}
\newcommand{\definition}[1]{\begin{defi}#1\end{defi}}
\newcommand{\coro}[1]{\begin{cor}#1\end{cor}}

\newcommand{\enumer}[1]{
                        \renewcommand{\labelenumi}{\arabic{enumi}.}
                        \begin{enumerate}#1\end{enumerate}
                        }

\newcommand{\feld}[2]{\begin{array}{#1} #2 \end{array}}
\newcommand{\matr}[2]{ \left( \begin{array}{#1} #2 \end{array} \right) }
\newcommand{\platz}{\vspace{0,1cm}}
\catcode`\"=\active
\catcode`\=\active
\catcode`\"=\active
\catcode`\Ž=\active
\catcode`\š=\active
\catcode`\™=\active
\catcode`\á=\active
\def"{\"a}
\def{\"u}
\def"{\"o}
\defŽ{\"A}
\defš{\"U}
\def™{\"O}
\letá=\ss

\begin{document}

\titel{spinor equations in weyl geometry\footnote{Supported by the SFB 288 of the DFG}}
\autor{Volker Buchholz\footnote{ e-mail: {\tt bv@mathematik.hu-berlin.de}}}
\adresse{ Humboldt Universit"at zu Berlin, Institut f"ur Reine Mathematik,\\ Ziegelstra"se 13a , D-10099 Berlin. } 
\abstrakt{In this paper, the Dirac, twistor and Killing equations on Weyl
manifolds with CSpin structures are
investigated. A conformal Schr"odinger-Lichnerowicz formula is presented and 
used to show integrability conditions for these equations. By introducing
 the Killing equation for spinors
of arbitrary weight, the result
of Andrei Moroianu in \cite{mor96} is generalized in the following sense. The
only non-closed Weyl manifolds of dimension greater than 3 that admit solutions of the real Killing equation are
4-dimensional and non-compact. Any Weyl manifold of these dimensions admitting a real Killing spinor
 has to be Einstein-Weyl.}{Differential Geometry}{53C05;53C10;53A30}
{ Weyl geometry; Dirac equation; twistor equation}

\section{Introduction}

In the first section, we state some basic definitions of density bundles, Weyl structures, curvature terms and
Einstein-Weyl structures.

The second section is dedicated to Dirac- and twistor operators on Weyl
manifolds.
In \cite{frie89}, \cite{bfgk}, \cite{lichn87} and  \cite{lichn89} 
the  properties and integrability conditions of the twistor and Killing equation
were intensively studied in the context of Riemannian geometry. Here we want to generalize some of these results to arbitrary Weyl
structures and spinor fields of arbitrary weight.
 The first result
in this area is due to Andrei Moroianu \cite{mor96} and deals with the
integrability conditions for the existence of non-trivial parallel spinors of
weight $0$. He found that the given Weyl structure has to be flat (closed)
on manifolds, which are not 4-dimensional and non-compact. Furthermore, he
gave several counter examples by  showing, that
in dimension 4 the existence of a parallel spinor field is equivalent to the
existence of a hypercomplex structure. This means in particular (see \cite{peswa}), that the Weyl structure
is Einstein-Weyl. We generalize this result to any dimension $n>2$ as well as
for Killing spinor fields $\psi$ of arbitrary weight satisfying
 $$\nabla^S_X\psi=\beta X\cdot\psi,$$ 
 where $\beta$ denotes a
complex density of weight $-1$. The
only non-closed Weyl manifolds of dimension greater than 3 that admit solutions of the real Killing equation are
4-dimensional and non-compact. Any Weyl manifold of these dimensions admitting a real Killing spinor
 has to be Einstein-Weyl. To this end it is crucial to
proof a generalized Schr"odinger-Lichnerowicz formula:
$$
 {\cal D}^2\psi=\Delta^S\psi+ \frac{1}{4}R\psi
+\left(\frac{n-2+2w}{4}\right)F\cdot\psi,
$$
where $R$ denotes the scalar curvature and $F$ the Faraday curvature of the Weyl
structure.   

This formula is also used in order to investigate integrability conditions of
the twistor equation, which is defined by  
$$
0=\twist_W\psi:=\nabla^S\psi+\frac{1}{n}\nu\dirac\psi.
$$ 
We then compute
$$
\nabla^S\dirac\psi=
\frac{n}{n-2}\left[-\frac{1}{2}\mu^2Ric'+\frac{1}{4(n-1)}R\nu
+\left(w-\frac{1}{2}\right)\left(\mu^2F+\frac{2}{4(n-1)}\nu
\mu F\right)\right]\psi
$$
on its kernel, where $Ric'$ is the Ricci curvature of the $\la{o}(n)$-component $W'$ of the Weyl structure W.
This equation corresponds to 
the equation $\nabla^S_X\dirac\psi=\frac{n}{2(n-2)}\left(\frac{R}{2(n-1)}X-Ric(X)\right)\cdot \psi$ in \cite{frie89}. We use this result
to prove, that the two well known first
integrals $C(\psi)$ and $Q(\psi)$ are parallel densities if the weight of 
$\psi$ is $\frac{1}{2}$ or $d \theta \cdot \psi = 0 $. Furthermore, we use it in
order show that the zeros of a twistor spinor field form a discrete set.

\vspace{0,5cm}

I would like to thank Thomas Friedrich for numerous discussions and hints on
this subject.

\section{Weyl geometry on conformal Spin manifolds}

Let $M^n$ be a smooth, oriented manifold and $({ \bf R}, M^n, \pi, GL(n,\r)$ its
frame bundle. Let $CO(n)_+=SO(n)\times \r_+$. For a conformal class $c$ let $ {\bf P}$ denote the corresponding
$CO(n)_+$-reduction.  
We define a two-fold covering $\map{\lambda^c}{Spin(n)\times \r_+=:CSpin(n)}{CO(n)_+=\{A\in
CO(n)|det(A)>0\}}$ by
$$
\lambda^c(a,\vartheta):=\vartheta \lambda(a),
$$
$\map{\lambda}{Spin(n)}{SO(n)}$ is the covering of the $SO(n)$.
The spinor representation $\kappa^w$ of $CSpin(n)$ on
$\Delta_n:=\com^{2^{\left[\frac{n}{2}\right]}}$ with
weight $w$ is defined as follows:
$\kappa^w(a,\vartheta)=\vartheta^{w}\kappa(a)$, where $\kappa$ is the
$Spin(n)$-representation on $\Delta_n$.
Like a Spin structure  a CSpin structure on $(M^n,c)$ is a pair
	$({ \bf P}_{CSpin},\Lambda^c)$, where \newline $({ \bf P}_{CSpin},\pi_{CSpin},M^n, CSpin(n))$ is a
	$CSpin(n)$-principal fibre bundle
   	on $M^n$ and \newline $\map{\Lambda^c}{{ \bf P}_{CSpin}}{{ \bf P}_{CO_+}}$ is a two-fold
   	covering that commutes with $\lambda^c$ and the action of the structure
	group. 
The existence of $CSpin$ structures is equivalent to the existence of $Spin$ structures, since $Spin(n)$ is maximally compact in $CSpin(n)$.   
We have the following vector bundles:
\enumer
      {
      \item ${\cal L}^w:={ \bf P}_{CSpin}\times_{|det\circ\lambda^c|^{\frac{w}{n}}}\r $ is called  density bundle  
             with weight $w$. 
      \item $(T^{r,s})^w:={ \bf
	P}_{CSpin}\times_{(\rho^{r,s}\circ\lambda^c)^w}\left(\bigotimes^r\dual\bigotimes^s\r^n\right)$,
	is the $(r,s)$-Tensor bundle with weight $w$. $(\rho^{r,s}\circ\lambda^c)^w$ denotes the
	the standard representation of $CSpin(n)$ on the $(r,s)$-tensors with
	weight w, 
	$(\rho^{r,s}\circ\lambda^c)^w(a,\vartheta)=\vartheta^w(\rho^{r,s}\circ\lambda^c)(a)$.
	$T$ shall denote the ordinary tagent bundle and $T^*$ its dual. 
      \item $S^w:={ \bf P}_{CSpin}\times_{\kappa^w}\Delta_n$ is the spinor
       bundle with  weight $w$. $S:=S^{1}$ denotes the ordinary spin bundle,
      }
Let $|vol_g|^{-\frac{1}{n}}=:l_g\in \dichte^1$ denote the density
corresponding to a metric $g \in c$. Then we are given the following conformally
invariant operators     
\enumer
      {
      \item $\map{c:=l_g^2g}{T^w\otimes T^{w_1}}{\dichte^{w+w_1}}$, $|X|^2:=
	c(X,X)\in
	\dichte^{2w}$ for $X\in T^w$. 
      \item $\map{(.)_c:=l_g^{2}(.)_g}{T^w}{(T^*)^w}$,
	\item $\map{tr:=l_g^{-2}tr_g}{(T^{r,s})^w}{(T^{r-2,s})^w},$ $r\geq 2$.	      
      \item the conformal, hermitian product $\map{(.,.):=l_g^2(.,.)_g}{\Gamma(S^w\otimes
   	S^{w_1})}{\Gamma(\dichte^{w+w_1})}$ and
      \item the conformal Clifford product
      $\map{\mu:=l_g\mu_g}{\Gamma(T^w \otimes S^{w_1})}{\Gamma(S^{w+w_1})}$,
      }
where $(.,.)_g$ and $\mu_g$ denote the hermitian product and the Clifford product given on $(M^n,g\in c)$. 
We can use $(.)_c$ in order to define the $\mu$ on arbitrary $(r,s)$-tensor fields.  
The operator $\map{\mu^{ab}}{T^{r,s} \otimes S^w}{S^w}$ is the conformal Clifford product of a
spinor field of weight $w$ with the $b^{th}$ and then with the $a^{th}$ component of a tensor field. 
Example:
$$
\mu^{21}\gamma\otimes X \otimes\omega\otimes \psi = \omega \otimes X \cdot \gamma \cdot \psi,
$$
where $\gamma \otimes X \in (T^{1,1})^w$, $\omega\in T^{2,0}$ and $\psi \in \Delta_n^{w_1}$. 
Whenever there are no indices, the Clifford product ranges over all
components of the corresponding tensor, i.e. $\mu
A\otimes\psi=A\cdot\psi:=\sum_{i_1,...,i_r}A(e_{i_1},...,e_{i_r})e_{i_1} \cdots
e_{i_r}\cdot\psi$. The
operator $\map{\nu}{\Delta_n^w}{T^{1,0}\otimes \Delta_n^w}$ is defined as follows:
$$
X\rfloor\nu \psi=\mu X\otimes \psi=X\cdot \psi.
$$
Then $\mu\nu=-n$ holds. Some well known identities have then the
following appearance:
\begin{eqnarray}
\mu^{12}\omega\otimes\psi & =
&-\mu^{21}\omega\otimes\psi-2tr^{12}\omega\psi\label{vertauschen}\\
tr\nu\omega\psi & = & w\cdot\psi\label{spurcliff}\\
Re(\nu \psi,\nu \psi) & = & (\psi,\psi)c:=|\psi|^2c\label{reellrausziehen}
\end{eqnarray}
Moreover, we define some operators on $T^{r,s}$:
\enumer
	{
	\item Let $(ab)$ denote the transposition of the components $a$ and $b$,
	\item $Sym:= Id + (12)$, $Alt = Id - (12)$, $Zyk:=Id +
	(23)(12)+(12)(23),$ $Zyk^{1234}:=Id+ (12)(23)(34)+(34)(23)(12)+(13)(24)$.
	}
A torsion-free connection $\map{W}{T{\bf P}}{\la{co}(n)}$ on a
   conformal manifold $(M^n,c)$ is called Weyl structure.
   $\nabla$ shall denote the induced covariant derivatives on
   associated vector bundles.   
The operators $c$, $tr$ and $(.)_c$ are parallel with respect to any
Weyl structure.
On $\dichte^1$ the curvature of a Weyl structure is given by $Alt\nabla^{T^*\otimes\dichte^1}\nabla^{\dichte^1}=:F \in \Omega^2(M)$. This
globally defined 2-form is called Faraday curvature.    
Choosing a gauge $g$ on $M^n$ provides a 1-form $\theta\in \Omega^1(M^n)$ in
the following way: $\nabla l_g=\theta\otimes l_g.$
For any gauge, we obtain $F=d\theta$.
Since the Lie algebra of the conformal
group splits into two components, there is also a splitting of a Weyl structure into 
a metric part $W'$ and a scalar part $\theta'\otimes Id$.
$$
W=W'+\theta'\otimes Id, \quad \map{W'}{T{\bf P}}{\la{o}(n)}, \quad
\theta'\in \Omega^1(P)
$$
A Weyl structure $W$ is exact (closed) if and only if 
$\theta$ is exact (closed) with respect to any gauge.

For a given Weyl manifold $(M^n,c,W)$ the curvature tensor   $\krumm \in
\Gamma(T^{4,0})^{-2}$ is defined
   by
   $$
   \krumm(X,Y,Z,U):=c(\nabla_X\nabla_YU
   -\nabla_Y\nabla_XZ-\nabla_{[X,Y]}Z,U),
   $$
   for any vector fields  $X,Y,Z,U \in\Gamma(T)$.  
The Ricci curvature is given by:      
$$
   Ric:=tr^{14}\R \mbox{ and } Ric':=tr^{14}\R'=Ric+F,
$$
where the primed objects belong to the connection $W'$.
$Ric$ is {\em not} symmetric. In fact, we obtain:
$\frac{1}{2}AltRic=-\frac{n}{2}F.$
Finally, we define the scalar curvature 
   $$
   R:=tr(Ric)=trtr^{14}\R \in \Gamma({\cal L}^{-2}),
   $$
which is not a function, but a density of weight -2.
\begin{lem}[Symmetry properties] Let $\R$' be the curvature tensor of $W'$.
Then
\glei{zyk1234}
{
\R' = (13)(24)\R'+\left[(13)+(23)-(14)-(24)\right]F\otimes c. 
}
\end{lem}
\beweis 
We have $Zyk\R'=-ZykF\otimes c$, which is just a version of the first Bianchi
identity for $\krumm$. This yields: 
$$ 
Zyk\R'Zyk^{1234}=-Zyk F\otimes cZyk^{1234}.
$$
We choose vector fields $X,Y,Z,T$. Then:
 \begin{eqnarray*}
  \lefteqn{Zyk\R'Zyk^{1234}(X,Y,Z,T)
=\R'(X,Y,Z,T)+\R'(Y,Z,X,T)+\R'(Z,X,Y,T)}\\
&+&\R'(Y,Z,T,X)+\R'(Z,T,Y,X)+\R'(T,Y,Z,X)
+\R'(Z,T,X,Y)+\R'(T,X,Z,Y)\\ &+&\R'(X,Z,T,Y)
+\R'(T,X,Y,Z)+\R'(X,Y,T,Z)+\R'(Y,T,X,Z)\\
 &=& 2\left(\R'(Z,X,Y,T)+\R'(T,Y,Z,X)\right)
=2(12)(23)(\R'-(13)(24)\R')(X,Y,Z,T).
\end{eqnarray*}
Similarly, we get for $F\otimes c$:
 $
  Zyk F\otimes cZyk^{1234}=2 Zyk^{1234} F\otimes c.
 $
Putting all this together implies:
 $$
 \R'=(13)(24)\R'-(12)(23)Zyk^{1234} F\otimes c
 =(13)(24)\R'+[(13)+(23)-(14)-(24)]F\otimes c,
 $$
 since 
 \glungen
   {
   (12)(23)Zyk^{1234} F\otimes c
   &=& (12)(23)[Id+ (12)(23)(34)+(34)(23)(12)+(13)(24)]F\otimes c\\   
   &=&[-(13)-(13)(23)(12)+(14)+(23)(24)(34)]F\otimes c\\
   &=& [-(13)-(23)+(14)+(24)]F\otimes c. 
   }
\schluss
\platz

For $n \geq 3$ a Weyl structure $W$ on $(M^n,c,W)$ is said to be Einstein-Weyl
if and only if
\glei{we}
        {
        Ric = \frac{R}{n} \cdot c - \frac{n}{2}F \quad or 
	\quad Ric'= \frac{R}{n} \cdot c - \frac{n-2}{2}F.
        }
The symmetric part of $Ric$ reduces to its trace if and only if $W$ is
Einstein-Weyl. $(M^n,c,W)$ is called an Einstein-Weyl manifold.

Let $W^S$ be the lift of $W$ into the  $CSpin$ structure and denote its induced
covariant derivative on $S^w$ by $\nabla^S$.
\begin{satz}{\em \cite{gau97}}\label{D genau} Fix a gauge  $g\in c$ and a spinor
$\psi\in \Gamma(S^w)$. Then the difference between
the spinor derivatives of $W$ and $W^g$, the Levi-Civita connection is as
follows:
\[
\nabla^{S}_X\psi-\nabla^{S,g}_X\psi=
 -\frac{1}{2}X\cdot\theta\cdot\psi+\left(w-\frac{1}{2}\right)\theta(X)\psi.
\]
\end{satz}
We define the spinorial curvature by
$
\krumm^{S,w}:=Alt\nabla^{T^*\otimes S}\circ\nabla^S=\kappa^w_*\Omega^S,
$
where $\Omega^S$ is the curvature
form of $W^S$. 
\begin{lem} \lab{spinkrumm}
\begin{eqnarray}\lab{R^S}  
{\cal R }^{S,w} & = &\frac{1}{4}\mu^{34}{\cal R }'  + w F   \\
\lab{mu234}  
\mu^{234}{\cal R }'&=&-2\mu^2Ric'-2\mu^2F-\nu\mu F \\ 
\lab{R'cliff}
\mu{\cal R }'&=&2R+2(n-2)\mu F  
\end{eqnarray}  
\end{lem}  
\beweis (1) According to the splitting $W=W'+\theta'\otimes Id$ we get for the corresponding curvature form:
$\Omega= \Omega' + F$. 
\begin{eqnarray*}
{\cal R }^S& = &\kappa_*^w  
\Omega^S=\kappa_*^w(\lambda^c_*)^{-1}\Omega  
=\kappa^w_*(\lambda_*^c)^{-1}\Omega'+wF
=\frac{1}{4}\mu^{34}{\cal R }'  + w F'.
\end{eqnarray*}  

(2) First we  use the symmetry properties of
$F\otimes c$ and (\ref{vertauschen}) to calculate:
\begin{eqnarray*}
\lefteqn{-\mu^{124}[(13)+(23)-(14)-(24)]F \otimes c}\\ 
& = & \left[-\mu^{324}-\mu^{134}+\mu^{124}(14)-\mu^{124}-2\mu^1tr^{24} \right]
F\otimes c\\
&=&\left[2\mu^{234}+2\mu^2tr^{23}-\mu^{214}(14)-2tr^{12}\mu^4(14)
-\mu^{124}+2\mu^2tr^{23} \right]
F\otimes c\\
&=&\left[2\mu^{234}+4\mu^2tr^{23}+\mu^{214}+2\mu^1tr^{14} -2\mu^1tr^{23}
-\mu^{124} \right]
F\otimes c\\
&=&\left[2\mu^{234}+4\mu^2tr^{23}-\mu^{124}+4\mu^2tr^{23} -\mu^{124} \right]
F\otimes c\\
&=&\left[2\mu^{234}+8\mu^2tr^{23}-2\mu^{124} \right]
F\otimes c
=(-2n+8)\mu^2F-2F\cdot\nu
\end{eqnarray*}
and
\begin{eqnarray*}
\lefteqn{-\mu^{234}ZykF \otimes c 
 = -\mu^{234}\left[Id+(23)(12)+(12)(23)\right]F\otimes c}\\
 &=&-\left[\mu^{234}-\mu^{234}(23)+\mu^{124}\right]F\otimes c
=-\left[2\mu^{234}+2\mu^2tr^{23}+\mu^{124}\right]F\otimes c\\
&=&(2n-2)\mu^2F-F\cdot\nu.
\end{eqnarray*}
This implies, by using $F\cdot\nu=\nu\mu F + 4 \mu^2F$ and equation (\ref{zyk1234}):
\begin{eqnarray*}  
\lefteqn{\mu^{234} \R' = -\mu^{234}\left[\left[(12)(23)+(23)(12)\right]
{\cal R }'-ZykF\otimes c\right]}\\
&=&\mu^{234}(23)\R' -\mu^{124}\R'-\mu^{234}ZykF\otimes c\\
&=&\mu^{234}(23)\R' -\mu^{124}(13)(24)\R'-\mu^{124}
[(13)+(23)-(14)-(24)]F\otimes c-\mu^{234}Zyk F\otimes c\\ 
&=&\mu^{234}(23)\R' +\mu^{324}\R'+2\mu^2tr^{24}\R'+(-2n+8)\mu^2F-2F\cdot\nu+
(2n-2)\mu^2F-F\cdot\nu\\  
&=&-2\mu^{234}\R' -6\mu^2Ric'+6\mu^2F-3F\cdot\nu
=-2\mu^{234}\R' -6\mu^2Ric'-6\mu^2F-3\nu\mu F.
\end{eqnarray*}
(3) From (\ref{mu234}), we obtain
\begin{eqnarray*}
\mu{\cal R }' &=&\mu\mu^{234}{\cal R
}'-2\mu Ric'-
\mu F+n\mu F
 = -\mu( Alt(Ric')+Sym(Ric'))+(n-2)\mu F\\
& = &(n-2)\mu F+2R+(n-2)\mu F
 = 2(n-2)\mu F+2R.
\end{eqnarray*}
\schluss  

\section{Spinor equations in Weyl geometry}   

\subsection{The Dirac operator}

\begin{defi}[Dirac operator]
The Dirac operator  
${\cal D}_W: \Gamma(S^w) \longrightarrow \Gamma(S^{w-1})$  is  defined
             by:  
             $$ {\cal D}_W:= \mu \nabla^{S}.$$ 	     
\end{defi}   
If there is no ambiguity to be expected, we will omit the index $W$.

\definition
   {
   The spinor Laplacian is given by
   $\map{\Delta^{S,w}}{\Gamma(S^w)}{\Gamma(S^{w-2})}$
   $$
   \Delta^{S,w}:=-tr\nabla^{T^*\otimes S}\circ\nabla^S.
   $$
   }

\begin{satz}[Schr"dinger-Lich\-nero\-wicz formula]  
Let  $\psi\in\Gamma(S^w)$. Then
\begin{equation}\lab{lichner}  
{\cal D}^2\psi=\Delta^S\psi+ \frac{1}{4}R\psi
+\left(\frac{n-2+2w}{4}\right)F\cdot\psi.
\end{equation}
\end{satz}  
\beweis This can be obtained directly. The final reduction of the curvature terms is due to the equations (\ref{R^S}) and
(\ref{R'cliff}).
\begin{eqnarray*} 
{\cal D}^2\psi &= & \mu\nabla^S\mu\nabla^S\psi
=\mu\nabla^{T^*\otimes S}\nabla^S\psi
=\mu \frac{1}{2}(Alt\nabla^{T^*\otimes S}\nabla^S+Sym \nabla^{T^*\otimes S}\nabla^S) \psi\\
& = &\mu \frac{1}{2}{\cal R }^S\psi  
- tr \nabla^{T^*\otimes S} \nabla^S \psi
=\Delta^S \psi+\frac{1}{8}\mu {\cal R' }\psi  
+ \frac{1}{2}w F \cdot \psi\\
&=&  \Delta^S
\psi+\frac{1}{4}R\psi+\left(\frac{n-2}{4} +
\frac{w}{2}\right)F\cdot\psi.
\end{eqnarray*}  
\schluss

\subsection{The twistor operator}

\begin{defi}[twistor operator]  
We define the twistor operator $\map{{\cal T}_W}{\Gamma(S^w)}{\Gamma(T^*
\otimes S^w)}$ of a $CSpin$  manifold $(M^n,c,W)$ by
$
\twist_W:= \nabla^{S,w} + \frac{1}{n}\nu{\cal D}.
$
\end{defi}

Let $(M^n,c,W)$ be a $CSpin$ manifold and $\psi\in \Gamma(S^w)$ a twistor
spinor field, i.e. an element of the kernel of $\twist$.
Then
$\nabla^S\psi=- \frac{1}{n}\nu\dirac\psi $ is true and therefore
\begin{equation}  
\Delta^S\psi=-tr\nabla^{T^*\otimes S}\nabla^S\psi=\frac{1}{n}tr\nabla^{T^*\otimes S}\nu\dirac\psi  
=\frac{1}{n}tr\nu
\nabla^S\dirac\psi=  
\frac{1}{n}\dirac^2\psi
\end{equation}
is satisfied.
From the Schr"odinger-Lichnerowicz formula we obtain:
\begin{equation}\lab{lichnertwist}  
\dirac^2= \frac{n}{4(n-1)}R+ \frac{(n-2+2w)n}{4(n-1)}\mu F.
\end{equation}
This leads  to
\begin{satz}
Let $\psi\in \Gamma(S^w)$ be a twistor spinor. Then:
\begin{equation} \lab{nabladirac} 
\nabla^S\dirac\psi=
\frac{n}{n-2}\left[-\frac{1}{2}\mu^2Ric'+\frac{1}{4(n-1)}R\nu
+\left(w-\frac{1}{2}\right)\left(\mu^2F+\frac{1}{2(n-1)}\nu
\mu F\right)\right]\psi.
\end{equation}
\end{satz}  
\beweis We use (\ref{mu234}) and (\ref{R^S}) in the first and second step. Then finally, after some direct calculations, we use
(\ref{lichnertwist}).   
\begin{eqnarray*}  
\lefteqn{- \frac{1}{2} \mu^2Ric'\psi -\frac{1}{2}\mu^2F\psi-\frac{1}{4}\nu\mu F\psi
= \frac{1}{4}\mu^{234}\R'}\\
&=&\mu^2{\cal R }^S\psi-w\mu^2F\psi
= \mu^2Alt \nabla^{T^*\otimes S}\nabla^S\psi-w\mu^2F\psi\\  
&=&- \frac{1}{n} \mu^2Alt \nabla^{T^*\otimes S}\nu\dirac\psi-w\mu^2F\psi
= \frac{1}{n} \mu^2Alt \nu\nabla^S\dirac\psi-w\mu^2F\psi\\  
&=& \frac{1}{n}\mu^2\nu \nabla^S\dirac\psi-\frac{1}{n} \mu^1\nu\nabla^S \dirac\psi-w\mu^2F\psi
=-\frac{1}{n} \nabla^S\dirac^2\psi- \frac{2}{n}\nabla^S\dirac\psi + \nabla^S
\dirac\psi-w\mu^2F\psi\\
&=&-\frac{n}{4n(n-1)}R\nu\psi -\frac{(n-2-2w)}{4(n-1)}\nu\mu F\psi +
\frac{n-2}{n}\nabla^S\dirac\psi -w\mu^2F\psi.
\end{eqnarray*}  
\vspace{-1cm}\schluss
\sa
   {
   If the term
   $
   \left(w-\frac{1}{2}\right)\left(\mu^2F+\frac{2}{4(n-1)}\nu
\mu F\right)\psi
   $
   reduces to a {\em single} Clifford product or even vanishes, e.g. if
   $w=\frac{1}{2}$ or $F \cdot\phi=0$ is satisfied, the sections
   $$
   C(\psi):=Re(\psi,\dirac\psi)\in\Gamma({\cal L}^{2w-1})
   $$
   and
   $$
   Q(\psi):= |\psi|^2|\dirac\psi|^2- trRe(\dirac\psi,\nu\psi)^2 \in
   \Gamma({\cal L}^{4w-2})
   $$
   are W-parallel.
   }
\beweis
In (\ref{nabladirac}) there are only single Clifford products left. Then (\ref{reellrausziehen}) yields
$$
\nabla C(\psi) = Re(\nabla^S\psi, \dirac\psi)+Re(\psi, \nabla^S\dirac \psi)
=Re(-\frac{1}{n}\nu\dirac\psi,\dirac\psi)=0
$$
and
\begin{eqnarray*}
\nabla Q(\psi) &=& 2Re (\nabla^S\psi,\psi) Re(\dirac\psi,\dirac\psi)+2 Re(\psi,\psi)Re(\nabla^S \dirac\psi,
\dirac\psi) \\
& &- 2tr^{23}Re(\nabla^S \dirac\psi, \nu\psi)Re( \dirac\psi, \nu\psi) 
+ \frac{2}{n}tr^{13}Re(\dirac\psi, \nu\nu\dirac\psi)Re(\dirac\psi,
\nu\psi)\\
&=& 0.
\end{eqnarray*}
\schluss

\subsubsection{The zeros of a twistor spinor field}

In this section we show that the zeros of a
twistor spinor field are a discrete set in $M^n$. Let $M^n$ be connected.
 We define 
$ E^w:=S^{w} \oplus S^{w-1}$ and regard the covariant derivative
 $\nabla^{E^w},$ which is characterized by
$$
\nabla^{E^w}=\left(\feld{cc}{
                       \nabla^{S,w} & \frac{1}{n}\nu\\
                       K^w          &  \nabla^{S,w-1}
                      }
               \right),	     
$$
where $\map{K^w = \nabla^{S,w-1}\dirac}{\Gamma(S^w)}{\Gamma(S^{w-1})}$.
\sa
   {
   For all  twistor spinor fields $\phi \in \Gamma(S^w)$
   $$
   \nabla^{E^w}\matr{c}{\phi\\ \dirac\phi} =0
   $$
   holds.
   Conversely, any $\nabla^{E^w}$-parallel section
   $ \matr{c}{\phi\\ \psi}\in \Gamma(E^w)$ yields:
   $$
   \twist_W\phi=0 \quad {\em and} \quad \psi=\dirac\phi.
   $$
  }
Since parallel sections on vector bundles over connected manifolds are uniquely
determined by their value in a single point, we obtain:
\coro
   {
   The dimension of the space of all twistor spinor fields of
   connected  $CSpin$ manifold is less than or equal to
   $2^{[\frac{n}{2}]+1}$.
   Furthermore, a  twistor spinor field $\psi$ on a connected
   CSpin manifold for that $\psi(m)=0$ and $\dirac\psi(m)=0$
   in a point $m\in M^n$ is trivial.
   }
\sa
   {
   The set $N_{\psi}:=\{\psi\in \Gamma(S^w):\mbox{ }
   \psi(m)=0 \mbox{ and } \twist_W\phi=0\}$ of a
   twistor spinor field $0\neq \psi \in \Gamma(S^w)$ is discrete in $M^n$.
   }
\beweis
 (\ref{nabladirac}) yields  $\nabla^S\dirac
\phi(m)= 0$ and for $g \in c$:
$$
0=2Re(\nabla\psi,\psi)(m)=\nabla|\psi|^2(m)=\nabla^g|\psi|^2(m),
$$
since $\nabla l^w=\nabla^g\l^w+w\theta\otimes l^w$.
$|\psi|^2$ is a density, i.e a section of $\dichte^{2w}$. Therefore,
we get for vector fields $X,Y$ on $M^n$:

$$ \nabla_X\nabla_Y|\psi|^2(m)=\nabla^g_X\nabla_Y|\psi|^2(m)
   +2w\theta(X)\nabla_Y|\psi|^2(m)
   =\nabla^g_X\nabla^g_Y|\psi|^2(m).
$$
If we choose $X$ and $Y$ to be W-parallel in $m$, we finally obtain by applying (\ref{reellrausziehen}) 
\gleifeldnon
   {
   \nabla_X\nabla_Y|\psi|^2(m)=2\nabla_X(\nabla_Y\psi,\psi)(m) =-\frac{2}{n}\nabla_X(Y\cdot\dirac\psi,\psi)(m)
   }
   {
   &=&\frac{2}{n^2}(Y\cdot\dirac\psi,X\cdot\dirac\psi)(m)
   =\frac{2}{n^2}c(X,Y)|\dirac\psi|^2(m).
   }
The combination of the latter two equations yields that
 $Hess_m(|\psi|^2)$ is not degenerated if $\dirac\psi(m)$ is not trivial.
Therefore, $m$ is an isolated point of $N(\psi)$. Otherwise, it follows
from the last corollary that $\psi$ must be trivial.

\schluss

\subsection{The Killing equation}

In this section, $M^n$ shall be connected.
\definition
   {
   [Killing spinor fields]
   A spinor field $\psi\in\Gamma(S^w)$ is called a 
    Killing spinor field if it satisfies the following
    differential equation:
   $$
   \nabla^S\psi=\beta\nu\psi, \quad \beta\in \Gamma(\com \otimes
   \dichte^{-1}),
   $$
   where $\beta$ is the  Killing density of $\psi$.
   }
A non-trivial Killing spinor field vanishes nowhere
on a connected manifold,
since it is parallel with
respect to the covariant derivative $\nabla^S-\beta\nu$.
It is obvious, that any Killing spinor field satisfies
the twistor equation and can be taken as an eigenspinor
of the Dirac operator with the 
eigen density $-n\beta$.
We now investigate the integrability conditions for the existence of non-trivial
Killing spinor fields. 
\begin{satz} \lab{killingsatz}
   Let $\psi\in\Gamma(S^w)$ be a Killing spinor field. 
\begin{enumerate}
        \item $\beta$ pureley imaginary, $w \neq \frac{n-2}{2}$: $(F\cdot\psi,\psi)=0$.
        \item $\beta$ real: $R = 4n(n-1)\beta^2$. 
        \begin{enumerate}
                \item $w \neq 0$: $W$ is exact and Einstein-Weyl
                \item w=0:  
                \begin{enumerate}
                        \item $\beta\neq 0$, $n \geq 4$: $W$ is exact and Einstein-Weyl
                        \item $\beta = 0$, $n>2$, ($4 \neq n $ or $M$ compact): $W$ is closed and Einstein Weyl
                        \item   $\beta = 0$, $n=4$, $M$ non-compact: $W$ is
								Einstein-Weyl
and $F$ is harmonic. 
                \end{enumerate} 
        \end{enumerate}         
\end{enumerate}
\end{satz}
\bemerkung For the latter  case ($n=4$ and $M^n$ non-compact) Moroianu gave in
      \cite{mor96} an example of a $CSpin$ manifold together with a non-closed 
      Weyl structure that carries non-trivial parallel spinor fields.\\

\beweis

\begin{enumerate}
      \item $\beta$ is purely imaginary: We have
      \glei{killinglemma}
      {
      R\psi + 2\left(\frac{n-2}{2}+w\right)F\cdot\psi  =  4(n-1)n\beta^2
      \psi-4(n-1)\nabla\beta\cdot\psi,
      }
      which itself follows from  (\ref{lichnertwist}):
        $$
        \frac{n}{4(n-1)}R\psi+\frac{(n-2+2w)n}{4(n-1)}F\cdot\psi=\dirac^2\psi
        =-n\mu\nabla^S(\beta\psi)
        =-n\nabla\beta\cdot\psi+n^2\beta^2\psi.
        $$
      The imaginary part of the product of (\ref{killinglemma}) with $\psi$ is
      as follows: 
      $$
      \left(w+\frac{n-2}{2}\right)(F\cdot\psi,\psi)=0,
      $$ 
      i.e., we have shown the assertion.
 \item $\beta$ is real: By multiplying (\ref{killinglemma}) with $\psi$
 	we see that $R=4n(n-1)\beta^2$ holds. 
 	
 	 $w\neq 0$: We obtain 
      $
      \nabla_X (\psi,\psi)=(\nabla^S_X\psi,\psi)+(\psi,\nabla^S_X\psi)=\beta(X\cdot
      \psi,\psi)-\beta(X\cdot\psi,\psi)=0,
      $
      by assumption.
      Hence, $W$ is exact. This means, that there is a metric $g$ of the conformal
       class $c$, 
        for that $W$ is the Levi-Civita connection and $W$ admits a Killingspinor. 
        Therefore, $(M^n,g)$ is Einstein, hence $(M^n,c,W)$ is Einstein-Weyl.
        
         $w=0$: By using the definition of Killing spinor fields we obtain:
      $
      \krumm^S\psi=Alt(\nabla\beta)\nu\psi+2\beta^2(\nu^{21}+c)\psi,
      $
      where $X,Y\rfloor\nu^{21}\psi=Y\cdot X\cdot \psi$.
      This yields:
      $$
      \mu^2\krumm^S\psi=-n\nabla\beta\otimes\psi-\nabla\beta\cdot\nu\psi-2\beta^2(n-1)\nu\psi.
      $$
      Together with (\ref{mu234}) we obtain:
      \glei{ric1}
        {
        \mu^2Ric'\psi=2n\nabla\beta\otimes\psi+2\nabla\beta\cdot\nu\psi+4\beta^2(n-1)\nu\psi-\mu^2F\psi-\frac{1}{2}\nu F\cdot\psi.
        }
       By (\ref{killinglemma}) and the assumption $R=4n(n-1)\beta^2$ holds and thus we obtain again from (\ref{killinglemma}):
      \glei{nunablabeta}
        {
        F\cdot\psi=-\frac{4(n-1)}{n-2}\nabla\beta\cdot\psi.
        }
      Inserting (\ref{nunablabeta}) into (\ref{ric1}) yields:
      \glei{ric2}
        {
        \mu^2Ric'\psi=2n\nabla\beta\otimes\psi+2\nabla\beta\cdot\nu\psi+\frac{R}{n}\nu\psi-\mu^2F\psi+\frac{2(n-1)}{n-2}\nu \nabla\beta 
        \cdot\psi.
        }
      The operator $\nabla\beta\cdot\nu$ consists of double Clifford products and scalar parts. We now rearrange (\ref{ric2}) accordingly.
	$$
        \mu^2Ric'\psi = 2\left(n-1-\frac{n-1}{n-2}\right)\nabla\beta\otimes\psi+\left(1-\frac{n-1}{n-2}\right)
        \mu^{12}Alt\nabla\beta\otimes c\psi +\frac{R}{n}\nu\psi-\mu^2 F\psi.
       $$
        If we multiply this equation by $\psi$, we see that $\nabla\beta$ must vanish for $n\neq 3$ since all the other terms  are purely       imaginary.
        Therfore, $W$ is exact and as before, $W$ is Einstein-Weyl. 
      \item The equations (\ref{R^S}) and (\ref{mu234}) 
      together with the assumption  yield 
      $$
      \mu^2Ric'\psi=-\mu^2F\psi-\frac{1}{2}\nu F\cdot\psi.
      $$ 
      Then (\ref{killinglemma}), $R=0$ and  the assumption impose
      $F\cdot\psi$ to vanish. 
      Since $\psi$ vanishes nowhere however, 
      $
      Ric'=-F.
      $
      Therefore, the symmetric part of $Ric'$ reduces to its trace (which is $0$),
       i.e., $W$ is Einstein-Weyl. Hence, Theorem 3.6 in \cite{peca} is
      applicable, which yields all remaining assertions. 
    \end{enumerate}
\schluss

\subsection{Two dimensional examples}

\enumer{
	\item {\sl Killing spinor fields of weight $\frac{1}{2}$}:\\
	We can find imaginary Killing spinor fields $\psi$ of weight $\frac{1}{2}$ on
	$(\r^2,[g],x_1dx^2)$.
	Because of Theorem \ref{D genau} they have to be a solution  of
	\begin{eqnarray*}
	X(\psi)=\frac{1}{2} X \cdot \theta\cdot \psi + \beta_g X \cdot \psi
	=X\cdot\left(\frac{1}{2}x_1\partial_2+\beta_g\right)\cdot\psi,
	\end{eqnarray*}
	 where $X\in \r^2$ and $\Gamma(\com \otimes \dichte^1) \ni \beta=\beta_gl_g$ hold true. 
	 $(\partial_1,\partial_2)$ are said to be the standard basis
	in $\r^2$. If one uses the following representation of the Clifford algebra
	$$
	   \partial_1 \longmapsto \left(\begin{array}{cc}i & 0 \\ 0 & -i\end{array}\right),\quad
	   \partial_2 \longmapsto \left(\begin{array}{cc}0 &i \\ i & 0\end{array}\right),
	$$
	one obtains
	$$
	X\cdot\matr{lr}
	   {
	   \beta_g & \frac{i}{2}x_1\\
	   \frac{i}{2}x_1 & \beta_g
	   }
	   \psi_0=0.
	$$
	We find a non-trivial kernel of the matrix for all $X$ if and only if
	$ \beta_g=\pm \frac{i}{2}x_1$. An element of this kernel must be
	of the form $\psi_0=\matr{r}{a\\\mp a}$ with $a\in \com$. We obtain, just as
	stated in Theorem \ref{killingsatz}:
	$$
	(F\psi_0,\psi_0)=x_1\matr{lr}{0&-1\\1&0}\matr{c}{a\\\mp
	a}\cdot\matr{r}{a\\ \mp a}=x_1\matr{r}{\pm a\\a}\cdot
	\matr{r}{a\\\mp a}=0.
	$$
	\item {\sl Parallel spinor fields of weight $0$:}\\
	On 
	$(\r^2,[g],x_1dx^2)$ we have to solve:
	$$
	   X(\psi)  =  \frac{1}{2}X\cdot\theta\cdot\psi+\frac{1}{2}\theta(X)\psi
	    =  \frac{1}{2}X_1x_1 \partial_1\cdot\partial_2\cdot\psi
	$$
	for $\psi\in \Gamma(S^0)$, where
	$X=\sum_{i=1}^2X_i\partial_i$.\\
	We use the following
	representation of the Clifford algebra 
	$$ 
	   \partial_1 \longmapsto \left(\begin{array}{cc}0 & i \\ i & 0\end{array}\right)
	   \quad
	   \partial_2 \longmapsto \left(\begin{array}{cc}0 &-1 \\ 1 & 0\end{array}\right).
	$$
	Therfore
	$$
	   \partial_1\cdot\partial_2 \longmapsto \left(\begin{array}{cc}i &0 \\ 0 & -i\end{array}\right),
	$$   
	and so we are given
	$$
	   X(\psi^+)  = \frac{i}{2}X_1x_1\psi^+,
	   \quad
	   X(\psi^-) = -\frac{i}{2}X_1x_1\psi^-,
	$$
	where $\psi^+$ and $\psi^-$ correspond to the splitting $\Delta_2=\Delta_2^+\oplus\Delta_2^-$.
	It is, however, not difficult to determine the solution of this system.
	$$
	   \psi^+(x):=\exp\left(\frac{i}{4}x_1^2\right)\psi_0^+,\quad
	   \psi^-(x):=\exp\left(-\frac{i}{4}x_1^2\right)\psi_0^-,
	$$
	where $\psi_0^{\pm} \in \com$.
	}
\selectlanguage{\english}

\end{document}